\documentclass[letterpaper, twocolumn]{IEEEtran}
\usepackage{multirow}
\usepackage{booktabs}
\usepackage{amsfonts}
\usepackage{amsmath,cases}
\usepackage{amssymb}
\usepackage{graphicx}
\usepackage{cite}
\usepackage{epsfig}
\usepackage{url}
\usepackage{color}
\usepackage{algorithm}
\usepackage{algorithmic}
\usepackage{epstopdf}
\usepackage{color}
\usepackage[tight]{subfigure}
\usepackage [english]{babel}
\usepackage [autostyle, english = american]{csquotes}
\MakeOuterQuote{"}
\usepackage{color}
\usepackage{amssymb}
\usepackage{bm}
\usepackage{caption}
\usepackage{float}
\usepackage{multirow}% htensor-to-tensor projection://ctan.org/pkg/multirow
\usepackage{hhline}% htensor-to-tensor projection://ctan.org/pkg/hhline
\usepackage{algorithm,algorithmic}
\usepackage{array}
\usepackage{caption}
\ifCLASSOPTIONcompsoc
\usepackage[caption=false,font=footnotesize,labelfont=sf,textfont=sf]{subfig}
\else
\usepackage[caption=false,font=footnotesize]{subfig}
\fi
\usepackage{fixltx2e}
\usepackage{amsthm}
\theoremstyle{plain}

\newcommand{\Tr}{{\sf Trace}}
\newcommand{\la}{\langle}
\newcommand{\ra}{\rangle}

\newcommand{\al}{\alpha}

\newcommand{\mc}{\mathcal}
\newcommand{\mb}{\mathbb}

\newcommand{\tb}{\textbf}
\newcommand{\nn}{\nonumber}
\newcommand{\bea}{\begin{eqnarray}}
\newcommand{\eea}{\end{eqnarray}}
\newcommand{\beq}{\begin{equation}}
\newcommand{\eeq}{\end{equation}}

\newcommand{\bet}{\beta}

\newtheorem{ex}{Example}\newcommand{\Ex}{\begin{ex}\rm}
	\newcommand{\eex}{\end{ex}}

\begin{document}
\title{Efficient tensor completion: Low-rank tensor train}
\author{Ho N. Phien$^1$, Hoang D. Tuan$^1$, Johann A. Bengua$^1$\thanks{$^1$Faculty of Engineering and Information Technology,	University of Technology Sydney, Ultimo, NSW 2007, Australia; Email:ngocphien.ho@uts.edu.au, tuan.hoang@uts.edu.au, johann.a.bengua@student.uts.edu.au.} and Minh N. Do$^2$\thanks{
$^2$Department of Electrical and Computer Engineering and the Coordinated Science Laboratory, University of Illinois at Urbana-Champaign, Urbana, IL 61801 USA; Email: minhdo@illinois.edu }}

\maketitle

\vspace*{-1.0cm}

\begin{abstract}
This paper proposes a novel formulation of the tensor completion problem to impute missing entries of data represented by tensors. The formulation is introduced in terms of tensor train (TT) rank which can effectively capture global information of tensors thanks to its construction by a well-balanced matricization scheme. Two algorithms are proposed to solve the corresponding tensor completion problem. The first one called simple low-rank tensor completion via tensor train (SiLRTC-TT) is intimately related to minimizing the TT nuclear norm. The second one is based on a multilinear matrix factorization model to approximate the TT rank of the tensor and called tensor completion by parallel matrix factorization via tensor train (TMac-TT). These algorithms are applied to complete both synthetic and real world data tensors. Simulation results of synthetic data show that the proposed algorithms are efficient in estimating missing entries for tensors with either low Tucker rank or TT rank while Tucker-based algorithms are only comparable in the case of low Tucker rank tensors. When applied to recover color images represented by ninth-order tensors augmented from third-order ones, the proposed algorithms outperforms the Tucker-based algorithms.
\end{abstract}

	\title{}	
	\author{
		\IEEEcompsocitemizethanks{\IEEEcompsocthanksitem The authors are from the Faculty of Engineering and Information Technology,	University of Technology Sydney, Ultimo, NSW 2007, Australia.\protect\\ E-mail:
			ngocphien.ho@uts.edu.au, tuan.hoang@uts.edu.au.\protect}\\}	
\section{Introduction}
Tensors are multi-dimensional arrays, known as higher-order generalizations of matrices and vectors \cite{Tamara_2009}. Tensors provide a natural way to represent multi-dimensional data objects whose entries are indexed by several continuous or discrete variables.  Employing tensors and their decompositions to process data objects has become increasingly popular since
\cite{Vasilescu_2003,Sun_2005,Franz_2009}. For instance, a color image is a third-order tensor defined by two indices for spatial variables and one index for color mode. A video comprised of color images is a fourth-order tensor with an additional index for a temporal variable. Residing  in extremely high-dimensional data spaces, the tensors in practical applications are nevertheless often of
\emph{low-rank} \cite{Tamara_2009}. Consequently, they can be effectively  projected to much smaller subspaces underlying their decompositions such as the CANDECOMP/PARAFAC (CP)\cite{Carroll_1970,Harshman_1970}, Tucker \cite{Tucker_1966} and tensor train (TT) \cite{Oseledets_2011} or matrix product state (MPS)\cite{Fannes_1992,Klumper_1993,PerezGarcia_2007}. 

Motivated by the success of low rank matrix completion (LRMC) \cite{Cai_2010,Ma_2009,Recht_2010}, much recent effort has been made to extend its concept to low rank tensor completion (LRTC). In fact, LRTC has been seen pivotal in computer vision and graphics, signal processing and machine learning \cite{Ji2013,Signoretto_2010,Signoretto_2011,Gandy_2011,Tomioka_2011,Tan_2014,Xu}. The common target is to recover  missing entries of a tensor from its partially observed entities \cite{Bertalmio_2000,Komodakis_2006,Korah_2007}. Despite its potential application, LRTC remains a grand challenge due to the fact that minimization of the tensor rank, defined as CP rank \cite{Tamara_2009}, is  an NP-hard problem. There have been some progress in formulating the LRTC via \emph{Tucker rank} \cite{Ji2013,Gandy_2011}. However, a crucial drawback of Tucker rank is that its components are ranks of matrices constructed based on an unbalanced matricization scheme (one mode versus the rest). Therefore, the upper bound of each individual rank is often small, and may not be suitable for describing global information of the tensor, especially strongly correlated tensors of high orders. In addition, the matrix rank minimization are only efficient when the matrix is more balanced. As the rank of a matrix is not more than $\min\{n,m\}$, where $m$ and $n$ are the number of rows and columns of the matrix, respectively, the high ratio $\max\{m,n\}/\min\{m,n\}$ would effectively rule out the need of matrix rank minimization.

In this paper we introduce a novel formulation of LRTC in the concept of \emph{TT rank}\cite{Oseledets_2011} which is different from the Tucker rank. Specifically, the TT rank is constituted by ranks of matrices formed by a well-balanced matricization scheme, i.e. matricize the tensor along one or a few modes. This gives rise to a huge advantage for representing the tensor rank in the sense that its components can have large upper bounds. Consequently, it provides a much better means to capture the global information in the tensor. We will also analyze how the rank of a matrix is closely related to the concept of von Neumann entropy in quantum information theory\cite{Nielsen_2009}. Subsequently, it is shown that the weakness of LRTC formulated by Tucker rank can be mitigated by utilizing the TT rank.

Two algorithms are proposed to approximately solve the proposed LRTC formulation, namely SiLRTC-TT and TMac-TT. The former is based on the SiLRTC \cite{Ji2013} which employs the block coordinate descent (BCD) for optimization and tools such as singular value thresholding from the matrix rank minimization problem \cite{Cai_2010,Ma_2009}. The SiLRTC-TT mainly solves the \emph{TT nuclear norm} minimization problem that is a convex surrogate for the new LRTC. Here, we define TT nuclear norm of a tensor as a sum of weighted nuclear norms of matrices formed by matricizing the tensor along one or a few modes. The latter adapted from its counterpart, i.e. TMac \cite{Xu}, is related to solving a weighted multilinear matrix factorization model. Although this model is non-convex, it can be solved effectively in the sense that no computationally expensive SVD is needed compared to the former.

The algorithms are applied to complete both synthetic and real world data described by tensors which are assumed to have low Tucker rank or TT rank. Empirical results simulated by proposed algorithms for synthetic data are much more efficient than the Tucker-based algorithms in the case of low TT rank and are comparable in the case of low Tucker rank tensors. When studying real world data such as color images, we introduce a tensor augmentation scheme called \emph{ket augmentation} (KA) to represent a lower-order tensor by a higher-order one without changing its number of entries. The KA scheme extended from the one introduced in \cite{Lattore} provides a perfect means to obtain higher-order tensor representation of visual data. We apply the proposed algorithms to complete a few color images represented by ninth-order tensors and results show that our algorithms outperforms the Tucker-based ones. Besides, we will show that our proposed algorithms recover augmented tensors more efficiently than those without applying augmentation scheme.

The rest of the paper is organized as follows. Section \ref{sec1} provides some notations and preliminaries of tensors. In Section \ref{sec2}, we firstly review the conventional formulation of LRTC and then introduce our new formulation in terms of TT rank. The algorithms are then proposed to solve the problem. We introduce the tensor augmentation scheme KA in Section \ref{sec3}. Section \ref{sec4} presents simulation results. Finally, we conclude our work in Section \ref{sec5}.

\section{Notations and preliminaries of tensors \label{sec1}}
We adopt some mathematical notations and preliminaries of tensors in \cite{Tamara_2009}. A \textit{tensor} is a multi-dimensional array and its \textit{order} (also known as way or mode) is the number of its dimensions. Scalars are zero-order tensors denoted by lowercase letters $(x, y, z,\ldots)$. Vectors and matrices are the first- and second-order tensors which are denoted by boldface lowercase letters (\textbf{x}, \textbf{y}, \textbf{z},\ldots) and capital letters $(X, Y, Z,\ldots)$, respectively. A higher-order tensor (the tensor of order three or above) are denoted by calligraphic letters $(\mc{X}, \mc{Y}, \mc{Z},\ldots)$.

An \textit{N}th-order tensor is denoted as $\mc{X}\in\mathbb{R}^{I_1\times I_2\times\cdots\times I_N}$ where $I_k$, $k=1,\ldots, N$ is the dimension corresponding to mode $k$. The elements of $\mc{X}$ are denoted as $x_{i_1\cdots i_k\cdots i_N}$, where $1\leq i_k\leq I_k$, $k=1,\ldots, N$.

A mode-$n$ fiber of a tensor $\mc{X}\in\mathbb{R}^{I_1\times I_2\times\cdots\times I_N}$ is a vector defined by fixing all indices but $i_n$ and denoted by \text{\bf x}$_{i_1\ldots i_{n-1}:i_{n+1}\ldots i_{N}}$.

Mode-$n$ matricization (also known as mode-$n$ unfolding or flattening) of a tensor $\mc{X}\in\mathbb{R}^{I_1\times I_2\times\cdots\times I_N}$ is the process of unfolding or reshaping the tensor into a matrix $X_{(n)}\in\mathbb{R}^{I_n\times (I_1\cdots I_{k-1}I_{k+1}\cdots I_N)}$ by rearranging the mode-$n$ fibers to be the columns of the resulting matrix.  Tensor element $(i_1,\ldots, i_{n-1},i_n,i_{n+1},\ldots, i_{N})$ maps to matrix element $(i_n,j)$ such that
\bea
j=1+\sum_{k=1,k\neq n}^{N}(i_k-1)J_k~~\text{with}~~J_k=\prod_{m=1, m\neq n}^{k-1}I_m.
\label{indexj}
\eea

The mode-$n$ product of a tensor $\mc{X}\in\mathbb{R}^{I_1\times I_2\times\cdots\times I_N}$ with a matrix $A\in\mathbb{R}^{J\times I_n}$ results into a new tensor of size $I_1\times\cdots\times I_{n-1}\times J\times I_{n+1}\times\cdots\times I_N$ which is denoted as $\mc{X}\times_n A$. Elementwise, it is described by
\bea
(\mc{X}\times_n A)_{i_1\cdots i_{n-1}ji_{n+1}\cdots i_N}=\sum_{i_n=1}^{I_n}x_{i_1\cdots i_n\cdots i_N}a_{ji_n}.
\eea

The Tucker decomposition (TD) is a form of higher-order principle component analysis \cite{Tamara_2009,Tucker_1966}. It is employed to decompose a tensor into a core tensor multiplied by a matrix along each mode. In general, for a given tensor $\mc{X}\in\mathbb{R}^{I_1\times I_2\times\cdots\times I_N}$, its TD is written as,
\bea
\mc{X} = \mc{G}\times_1A^{(1)}\times_2A^{(2)}\cdots\times_NA^{(N)},
\eea
where the core tensor $\mc{G}\in\mathbb{R}^{r_1\times r_2\cdots\times r_{N}}$ and the factor matrices $A^{(k)}\in\mathbb{R}^{r_k\times I_{k}},k=1,\ldots,N$. The vector $\tb{r} = (r_1, r_2,\ldots, r_N)$, where $r_n$ is the rank of the corresponding matrix $X_{(n)}$ denoted as $r_n = \text{rank}(X_{(n)})$, is called as the \emph{Tucker rank} of the tensor $\mc{X}$.

The inner product of two tensors $\mc{X},\mc{Y}\in\mathbb{R}^{I_1\times I_2\times\cdots\times I_N}$ is defined as
\bea
\la\mc{X},\mc{Y}\ra &=&\sum_{i_1}\sum_{i_2}\cdots\sum_{i_N}x_{i_1i_2\cdots i_{N}}y_{i_1i_2\cdots i_{N}}.
\eea
Accordingly, the Frobenius norm of $\mc{X}$ is $||\mc{X}||_F = \sqrt{\la\mc{X},\mc{X}\ra}$.

\section{The formulation of tensor completion and algorithms \label{sec2}}
This section firstly revisits the conventional formulation of LRTC based on the Tucker rank, a generalization of LRMC. Then a new LRTC formulated in terms of the TT rank is introduced with algorithms.
\subsection{Conventional formulation of tensor completion}
We give an overview of matrix completion before introducing the formulation for tensor completion. The problem on how to recover missing entries of a low-rank matrix  $T\in\mb{R}^{m\times n}$ from its partially known entries given by a subset $\Omega$ can be studied via the well-known optimization problem \cite{Kurucz_2007}:
\bea
\begin{aligned}
	& \underset{X}{\text{min}}&&\text{rank}(X)~~~~ \text{s.t.}&& X_{\Omega} = T_{\Omega}.
\end{aligned}
\label{eq1}
\eea
The missing entries of $X$ are estimated such that the rank of $X$ is as small as possible. Due to the combinational nature of the function $\text{rank}(\cdot)$, problem (\ref{eq1}), however, is  NP-hard and one needs to look for its surrogates. Minimization of the
matrix nuclear norm has proved as an excellent approximation of the matrix rank. This leads to the following convex optimization problem for matrix completion \cite{Bach_2008,Cai_2010,Ma_2009}:
\bea
\begin{aligned}
	& \underset{X}{\text{min}}&&||X||_{*}~~~~ \text{s.t.}&& X_{\Omega} = T_{\Omega},
\end{aligned}
\label{eq2}
\eea
where the nuclear norm $||X||_{*}$ is the summation of the singular values of $X$. Note that this matrix nuclear norm optimization problem is efficient only when $X$ is balanced which implies that $m\approx n$.

Alternatively, one can also apply the low-rank matrix factorization model to solve the matrix completion problem\cite{Wen_2012}:
\bea
\begin{aligned}
	& \underset{U,V,X}{\text{min}}&&\frac{1}{2}||UV-X||^{2}_{F}~~~~ \text{s.t.}&& X_{\Omega} = T_{\Omega},
\end{aligned}
\label{eq3}
\eea
where $U\in\mb{R}^{m\times r}$, $V\in\mb{R}^{r\times n}$ and $X\in\mb{R}^{m\times n}$ and the integer $r$ is the estimated rank of the matrix $T$.

The matrix completion can be generalized for tensor underlying the concept of tensor rank, e.g. the CANDECOMP/PARAFAC rank (CP-rank), Tucker rank \cite{Tamara_2009,Tucker_1966}. For instance, in terms of Tucker rank, completing an $N$th-order tensor $\mc{T}\in\mathbb{R}^{I_1\times I_2\cdots\times I_{N}}$ from its known entries given by an index set $\Omega$ is related to solving the following optimization problem \cite{Ji2013,Gandy_2011,Tan_2014,Xu}:
\bea
\begin{aligned}
	& \underset{X_{(k)}}{\text{min}}&&\sum_{k=1}^{N}\al_{k}\text{rank}(X_{(k)})~~~~ \text{s.t.}&& \mc{X}_{\Omega} = \mc{T}_{\Omega}.
\end{aligned}
\label{eq4}
\eea
where $\{\al_k\}_{k=1}^{N}$ are defined as weights fulfilling condition $\sum_{k=1}^{N}\al_k=1$. The Eq.~(\ref{eq4}) is a weighted multilinear matrix completion problem which is still NP-hard. Therefore, one needs to switch to an alternative by generalizing the matrix case. For instance, (\ref{eq4}) can be converted to the following optimization problem \cite{Ji2013}:
\bea
\begin{aligned}
	& \underset{X_{(k)}}{\text{min}}&&\sum_{k=1}^{N}\al_{k}||X_{(k)}||_{*}~~~~ \text{s.t.}&& \mc{X}_{\Omega} = \mc{T}_{\Omega},
\end{aligned}
\label{eq4_1}
\eea
where $\sum_{k=1}^{N}\al_{k}||X_{(k)}||_{*}$ can be defined as \emph{Tucker nuclear norm} of the tensor. This problem can be then solved by simply applying known methods such as the block coordinate descent (BCD) to alternatively optimize a group of variables while the other groups remain fixed. Besides, one can generalize the low-rank matrix factorization model in (\ref{eq3}) as a subsitute for (\ref{eq4}) and then apply the BCD method with the nonlinear Gauss-Seidal method to solve it \cite{Wen_2012,Xu,Tan_2014}.

Although the Tucker-based LRTC problem has become increasingly popular, it is only appropriate for the tensors with either low Tucker rank or low orders and might be less efficient when applying to real world data represented by tensors of orders higher than three. This weakness comes from the fact that each matrix $X_{(k)}$ in (\ref{eq4}) is obtained by matricizing the tensor along one single mode. As a consequence, it is unbalanced and the corresponding rank is not large enough to capture the \emph{global correlation} between elements in the tensor. Even when all the modes have the same dimension ($I_1=\cdots=I_N\equiv I$), these matrices are highly unbalanced. We can clarify this observation via the concept of von Neumann entropy \cite{Nielsen_2009} as follows.

Represent $\mc{X}$ as a \emph{pure state} in the space $\mathbb{R}^{I_1\times I_2\cdots\times I_{N}}$,
\bea
\mc{X} = \sum_{i_1,i_2\ldots,i_N}x_{i_1i_2\cdots i_N}\tb{e}_{i_1}\otimes\tb{e}_{i_2}\cdots\otimes\tb{e}_{i_N},
\label{eq_state}
\eea
where "$\otimes$" denotes a tensor product \cite{Tamara_2009},  $\tb{e}_{i_k}\in\mathbb{R}^{I_k}$ forms an orthonormal
basis in $\mathbb{R}^{I_k}$ for each $k=1,\ldots,N$. Applying mode-$k$ matricization of $\mc{X}$ results into $X_{(k)}$ representing a \emph{pure state} of the composite system $AB$ in the space $\mc{H}_{AB}\in\mb{R}^{m\times n}$,
which is a tensor product of two subspaces $\mc{H}_{A}\in\mb{R}^{m}$ and $\mc{H}_{B}\in\mb{R}^{n}$ of dimensions $m = I_k$ and $n=\prod\limits_{l=1,l\neq k}^{N}I_l$, respectively. The subsystems $A$ and $B$ are
seen as two \emph{contigous partitions} consisting of mode $k$ and all other  modes of the tensor, respectively.
It follows from (\ref{eq_state}) that
\bea
X_{(k)} = \sum_{i_k,j}x_{i_kj}\tb{e}_{i_k}\otimes\tb{e}_{j},
\label{eq_state1}
\eea
where the new index $j$ is defined as in (\ref{indexj}), $\tb{e}_{j} = \otimes_{l=1,l\neq k}^{N}\tb{e}_{i_l}\in\mb{R}^{n}$. According to the Schmidt decomposition \cite{Nielsen_2009}, there exist orthonormal bases $\{\tb{u}^{A}_{l}\}$ in $\mc{H}_{A}$ and $\{\tb{v}^{B}_{l}\}$ in $\mc{H}_{B}$ such that,
\bea
X_{(k)} = \sum_{l=1}^{r_k}\lambda_l\tb{u}^{A}_{l}\otimes\tb{v}^{B}_{l},
\eea
where $r_k$ is the rank of $X_{(k)}$, $\lambda_l$ are nonvanishing singular values and $\{\tb{u}^{A}_{l}\}$ and $\{\tb{v}^{B}_{l}\}$ are columns of orthonormal matrices $U$ and $V$ which are obtained from the SVD $X_{(k)} = U\lambda V^T$, respectively. The correlation between two subsystems $A$ and $B$ can be studied via von Neumann entropy defined as \cite{Nielsen_2009}:
\bea
S^A=-\Tr(\rho^A\log_2(\rho^A)),
\label{en1}
\eea
where $\rho^A$ is called the \emph{reduced density matrix operator} of the composite system and computed by taking the partial trace of the density matrix $\rho^{AB}$ with respect to $B$. Specifically, we have
\bea
\rho^{AB} &=& X_{(k)}\otimes(X_{(k)})^T\nn\\
&=&\Big(\sum_{l=1}^{r_k}\lambda_l\tb{u}^{A}_{l}\otimes\tb{v}^{B}_{l}\Big)\otimes\Big(\sum_{j=1}^{r_k}\lambda_j\tb{u}^{A}_{j}\otimes\tb{v}^{B}_{j}\Big)^T.
\eea
Then $\rho^A$ is computed as
\bea
\rho^A &=& \Tr_B(\rho^{AB})\nn\\
&=&\sum_{l=1}^{r}\lambda^{2}_l\tb{u}^{A}_{l}\otimes(\tb{u}^{A}_{l})^{T},
\label{en2}
\eea
Substituting (\ref{en2}) to (\ref{en1}) yields
\bea
S^A= -\sum_{l=1}^{r_k}\lambda^{2}_l\log_{2}\lambda^{2}_l.
\eea
Similarly,
\bea
S^B&=&-\Tr(\rho^B\log_2(\rho^B))\nn\\
&=&-\sum_{l=1}^{r_k}\lambda^{2}_l\log_{2}\lambda^{2}_l,
\eea
which is the same with $S^A$, simply $S^A=S^B=S$. This entropy reflects the correlation or \emph{degree of entanglement} between  subsystem $A$ and its complement $B$ \cite{Bennett_1996}. Without loss of generality,  the normalization condition $\sum_{l=1}^{r_k}\lambda^{2}_l=1$ can be imposed, so $0\leq S\leq\log_{2} r_k$. Obviously, there is no correlation between subsystems
$A$ and $B$ whenever $S=0$ (where $\lambda_1 = 1$ and the other singular values are zeros).  There exists correlation between
 subsystems $A$ and $B$ whenever $S\neq 0$  with its maxima  $S=\log_{2} r_k$ (when $\lambda_1 = \cdots = \lambda_{r_k} = 1/\sqrt{r_k}$). If the singular values decay significantly, e.g. exponential decay, we can also keep a few $r_k$ ($r_k\ll m$) largest singular values of $\lambda$ without considerably losing accuracy in quantifying the amount of correlation between the subsystems. Then $r_k$ is referred to as the approximate low rank of the matrix $X_{(k)}$ which means that the amount of correlation between the elements in the matrix is small. On the contrary, if two subsystems $A$ and $B$ are highly correlated, i.e. the singular values decay very slowly, then $r_k$ needs to be as large as possible to capture the correlation in the tensor. Therefore, the problem of matrix rank minimization is in fact intimately related to the problem of von Neumann entropy minimization.

From the above analysis, we see that the amount of correlation between elements in the matrix $X_{(k)}$ depends on the rank $r_k$ which is bounded by $m=I_k$. Therefore, when the dimensions of modes are slightly different or the same, that is $I_1\approx I_2\approx\cdots\approx I_N\approx I$, the matrix $X_{(k)}$ is essentially unbalanced due to $m\ll n$ when either $I$ or $N$ is large. As a result, the limit of each $r_k$ is too small to describe the correlation of the tensor in case the tensor $\mc{X}$ has higher order ($N>3$) that makes the Tucker-based LRTC no longer appropriate for a highly-correlated tensor. In the next section we introduce a new LRTC problem formulated in terms of \emph{TT rank} defined by more balanced matrices.

\subsection{Tensor completion formulation in the concept of tensor train rank and algorithms}
The tensor train (TT) decomposition is applied to decompose a higher-order tensor into a sequence of connected lower-order tensors \cite{Oseledets_2011}. Using Vidal's decomposition \cite{Vidal_2004}, the TT decomposition of a tensor described by (\ref{eq_state}) can be written in the following form,
\bea
\mc{X}
&=&\sum_{i_1,\ldots,i_N}\Gamma^{[1]}_{i_{1}}\lambda^{[1]}\cdots\lambda^{[N-1]}\Gamma^{[N]}_{i_{N}}\tb{e}_{i_1}\otimes\cdots\otimes\tb{e}_{i_N},
\label{tt2}
\eea
where for $k=1,\ldots,N$, $\Gamma^{[k]}_{i_k}$ is an $r_{k-1}\times r_k$ matrix and $\lambda^{[k]}$ is the $r_{k}\times r_{k}$ diagonal singular matrix, $r_0=r_{N+1}=1$. For every $k$, the following orthogonal conditions are fulfilled:
\bea
\sum_{i_k=1}^{I_k}\Gamma^{[k]}_{i_{k}}\lambda^{[k]}(\Gamma^{[k]}_{i_{k}}\lambda^{[k]})^{T} &=& \mb{I}^{[k-1]},\\
\sum_{i_k=1}^{I_k}(\lambda^{[k-1]}\Gamma^{[k]}_{i_{k}})^{T}\lambda^{[k-1]}\Gamma^{[k]}_{i_{k}} &=& \mb{I}^{[k]},
\eea
where $\mb{I}^{[k-1]}$ and $\mb{I}^{[k]}$ are the identity matrices of sizes $r_{k-1}\times r_{k-1}$ and $r_{k}\times r_{k}$, respectively. Based on the form (\ref{tt2}), each component $r_k$ of the so-called \emph{TT rank} of the tensor, simply defined as $\tb{r} = (r_1,r_2,\ldots, r_{N-1})$, can be determined directly via the singular matrices $\lambda^{[k]}$. Specifically, to  determine $r_k$, rewrite (\ref{tt2}) as
\bea
\mc{X}
&=&\sum_{i_1,i_2\ldots,i_N}\tb{u}^{[1\cdots k]i_1\cdots i_k}\lambda^{[k]}\tb{v}^{[k+1\cdots N]i_{k+1}\cdots i_N},
\label{tt3}
\eea
where
\bea
\tb{u}^{[1\cdots k]i_1\cdots i_k}&=&\Gamma^{[1]}_{i_{1}}\lambda^{[1]}\cdots\Gamma^{[k]}_{i_{k}}\otimes_{l=1}^{k}\tb{e}_{i_l},
\eea
and
\bea
\tb{v}^{[k+1\cdots N]i_{k+1}\cdots i_N}&=&\Gamma^{[k+1]}_{i_{k+1}}\lambda^{[k+1]}\cdots\Gamma^{[N]}_{i_{N}}\otimes_{l=k+1}^{N}\tb{e}_{i_l}.~~~~~~
\eea
We can also rewrite (\ref{tt3}) in terms of the matrix form of an SVD as
\bea
X_{[k]}&=& U\lambda^{[k]}V^{T},
\label{tt4}
\eea
where $X_{[k]}\in\mb{R}^{m\times n}$ ($m = \prod_{l=1}^{k}I_l, n=\prod_{l=k+1}^{N}I_l$) is  the \emph{mode-$(1,2,\ldots, k)$ matricization} of the tensor $\mc{X}$ \cite{Oseledets_2011}, $U\in\mb{R}^{m\times r_k}$ and $V\in\mb{R}^{n\times r_k}$ are orthogonal matrices. Obviously, $r_k$, defined as number of nonvanishing singular values of $\lambda^{[k]}$, is the rank of $X_{[k]}$.

In practice, the mode-$(1,2,\ldots, k)$ matricization $X_{[k]}$ of tensor $\mc{X}$ can be obtained by reshaping the tensor $\mc{X}$ in such a way that the first $k$ indices enumerate the rows of $X_{[k]}$, and the last $(N-k)$ enumerate the columns of $X_{[k]}$ \cite{Oseledets_2011,Mu_2014}. Specifically, the tensor element $(i_1,i_2,\ldots, i_{N})$ of $\mc{X}$ maps to the element $(i,j)$ of
$X_{[k]}$ for
\bea
i &=& 1+ \sum_{m=1}^{k}\big((i_m-1)\prod_{l=1}^{m-1}I_l\big),\\
j &=& 1+ \sum_{m=k+1}^{N}\big((i_m-1)\prod_{l=k+1}^{m-1}I_l\big).
\eea
Since matrix $X_{[k]}$ is obtained by matricizing along a few $k$ modes rather than one single mode, its rank $r_k$ is bounded by $\min(\prod_{l=1}^{k}I_l,\prod_{l=k+1}^{N}I_l)$. Therefore TT rank is in general more appropriate than Tucker rank for quantifying correlation of  higher-order tensors.

We now propose to formulate the LRTC problem in terms of TT rank as
\bea
\begin{aligned}
	& \underset{X_{[k]}}{\text{min}}&&\sum_{k=1}^{N-1}\al_{k}\text{rank}(X_{[k]})~~~~ \text{s.t.}&& \mc{X}_{\Omega} = \mc{T}_{\Omega},
\end{aligned}
\label{eq5}
\eea
where $\al_{k}$ denotes the weight that the rank of the matrix $X_{[k]}$ contributes to the TT rank that the condition $\sum_{k=1}^{N-1}\al_k=1$ is satisfied. The LRTC problem is now relaxed to the weighted multilinear matrix completion problem which is similar to (\ref{eq4}). It is still difficult to directly tackle as $\text{rank}(\cdot)$ is presumably hard. Thus, we will convert this problem into two separate problems. The first one based on the so-called \emph{TT nuclear norm}, defined as
\bea
||\mc{X}||_{*}&=&\sum_{k=1}^{N-1}\al_{k}||X_{[k]}||_{*},
\label{TTnorm}
\eea
is given by
\bea
\begin{aligned}
	& \underset{\mc{X}}{\text{min}}&&\sum_{k=1}^{N-1}\al_{k}||X_{[k]}||_{*}~~~~ \text{s.t.}&& \mc{X}_{\Omega} = \mc{T}_{\Omega},
\end{aligned}
\label{eq6}
\eea
The problem (\ref{eq6}) is defined similarly to (\ref{eq4_1}) where the Tucker nuclear norm is used instead. Besides, from (\ref{eq6}) we can recover the square model \cite{Mu_2014} by choosing the weights such that $\al_{k}=1$ if $k=\text{round}(N/2)$ otherwise $\al_{k}=0$. 

The problem (\ref{eq6}) can be further converted to the following problem:
\bea
\begin{aligned}
	& \underset{\mc{X},M_k}{\text{min}}&&\sum_{k=1}^{N-1}\al_{k}||M_k||_{*} + \frac{\beta_k}{2}||X_{[k]}-M_k||^{2}_{F}\\
	& \text{s.t.}&& \mc{X}_{\Omega} = \mc{T}_{\Omega},
\end{aligned}
\label{eq7}
\eea
where $\bet_k$ are positive numbers and can be solved by employing the BCD method for the optimization which will be discussed later in this section.

The second proposed problem is based on the multilinear matrix factorization model. More specifically, given a  matrix $X_{[k]}\in\mb{R}^{m\times n}$ of rank $r_k$, it can be factorized as $X_{[k]} = UV$ where $U\in\mb{R}^{m\times r_k}$ and $V\in\mb{R}^{r_k\times n}$. Therefore, instead of optimizing the nuclear norm of the unfolding matrices $X_{[k]}$, the Frobenius norm is minimized as follows,

\bea
\begin{aligned}
	& \underset{U_k,V_k,\mc{X}}{\text{min}}&&\sum_{k=1}^{N-1}\frac{\al_{k}}{2}||U_kV_k-X_{[k]}||^{2}_{F} \\
	& \text{s.t.}&& \mc{X}_{\Omega} = \mc{T}_{\Omega},
\end{aligned}
\label{eq8}
\eea
where $U_{k}\in\mathbb{R}^{\prod_{j=1}^{k}I_{j}\times r_k}$ and $V_{k}\in\mathbb{R}^{r_k\times\prod_{j=k+1}^{N}I_{j}}$. This model is similar to the one proposed in \cite{Tan_2014,Xu} (which is an extension of the matrix completion model \cite{Wen_2012}) where the Tucker rank is employed.

To solve the convex but nondifferentiable optimization problem described by (\ref{eq7}), one can adapt the TT nuclear norm to the algorithms such as SiLRTC, FaLRTC in \cite{Ji2013}. Besides, in order to solve (\ref{eq8}), we can apply the alternating least squares (ALS) technique to variationally optimize $U,V,\mc{X}$ until a convergence is obtained. Specifically, one can modify the algorithms such as TMac and TC-MLFM in \cite{Xu} and \cite{Tan_2014}, respectively by incorporating the concept of TT rank into them. The essential advantage of this multilinear matrix factorization model when compared to the model in (\ref{eq7}) is that it avoids a lot of SVDs and hence can substantially save the computational time.

Let us propose the first algorithm to solve the optimization problem in (\ref{eq7}) which is deeply rooted by the SiLRTC algorithm \cite{Ji2013}.  We call our algorithm "SiLRTC-TT" which stands for "simple low rank tensor completion via tensor train". The central concept of this algorithm is based on the BCD method to alternatively optimize a group of variables while the other groups remain fixed. More specifically, the variables are divided into two main groups. The first one contains the unfolding matrices $M_1, M_2,\ldots, M_{N-1}$ and the other is tensor $\mc{X}$. Computing each matrix $M_k$ is related to solving the following optimization problem:
\bea
\begin{aligned}
	& \underset{M_k}{\text{min}}&&\al_{k}||M_k||_{*} + \frac{\beta_k}{2}||X_{[k]}-M_k||^{2}_{F},
\end{aligned}
\label{minrankX7}
\eea
with fixed $X_{[k]}$. The optimal solution for this problem has the closed form \cite{Ma_2009} which is determined by
\bea
M_k= \mathbf{D}_{\gamma_k}(X_{[k]}),
\eea
where $\gamma_k=\frac{\al_k}{\bet_k}$ and $ \mathbf{D}_{\gamma_k}(X_{[k]})$ denotes the thresholding SVD of $X_{[k]}$ \cite{Cai_2010}. Specifically, if the SVD of $X_{[k]} = U\lambda V^T$, its thresholding SVD is defined as:
\bea
\mathbf{D}_{\gamma_k}(X_{[k]}) = U\lambda_{\gamma_k}V^T,
\eea
where $\lambda_{\gamma_k} = diag(\max(\lambda_l-\gamma_k,0))$. After updating all the $M_k$ matrices, we turn into another block to compute the tensor $\mc{X}$ which elements are given by
\bea
x_{i_1\cdots i_N}= \left\{\begin{array}{ll}\Big (\frac{\sum_{k=1}^{N} \bet_k \text{fold}(M_k)}{\sum_{k=1}^{N}\bet_k}\Big)_{i_1\cdots i_N}&({i_1\cdots i_N})\notin\Omega\\
t_{i_1\cdots i_N}&({i_1\cdots i_N})\in\Omega\\
\end{array}\right.
\label{UpdateX}
\eea
The pseudo-code of this algorithm is given in Algorithm \ref{Algorithm 1}. The convergence condition is reached when the relative error between two successive tensors $\mc{X}$ is smaller than a threshold. The algorithm is guaranteed to be converged and gives rise to a global solution since the objective in (\ref{eq7}) is a convex and the nonsmooth term is separable. We can also apply this algorithm for the square model \cite{Mu_2014} by simply choosing the weights such that $\al_{k}=1$ if $k=\text{round}(N/2)$ otherwise $\al_{k}=0$. For this particular case, let us call the algorithm as SiLRTC-Square.
\begin{table}[!thb]
	\centering
	\caption{SiLRTC-TT}
	\label{Algorithm 1}	
	\begin{tabular}{*2l} % centered columns (4 columns)
		\hline
		~&~\\
		\tb{Input:} The observed data $\mc{T}\in\mathbb{R}^{I_1\times I_2\cdots\times I_{N}}$, index set $\Omega$.&\\
		\tb{Parameters:} $\al_k,\bet_k, k=1,\ldots,N-1$.&\\							
		~&~\\	
		\hline		
		~&~\\	
		1:~~\tb{Initialization:} $\mc{X}^{0}$, with $\mc{X}^{0}_{\Omega} = \mc{T}_{\Omega}$, $l=0$.& \\
		2:~~\tb{While not converged do:}&\\
		3:~~~~\tb{for} $k = 1$ \tb{to} $N-1$ \tb{do}&\\
		4:~~~~~~~Unfold the tensor $\mc{X}^{l}$ to get $X^{l}_{[k]}$&\\
		5:~~~~~~~$M^{l+1}_k = \mathbf{D}_{\frac{\al_k}{\bet_k}}(X^{l}_{[k]})$&\\	
		6:~~~~\tb{end for}&\\	
		7:~~ Update $\mc{X}^{l+1}$ from $M^{l+1}_k$ by (\ref{UpdateX})&\\
		8:~~\tb{End while}&\\
		\hline
		~&~\\
		\tb{Output:} The recovered tensor $\mc{X}$ as an approximation of $\mc{T}$&
	\end{tabular}
\end{table}

In order to solve the problem given by (\ref{eq8}), we apply the BCD method to alternatively optimize different groups of variables. Specifically, we can first solve the following problem:
\bea
\begin{aligned}
	& \underset{U_k,V_k,X_{[k]}}{\text{min}}&&||U_kV_k-X_{[k]}||^{2}_{F},
\end{aligned}
\label{minrankX7}
\eea
for $k=1,2,\ldots,N-1$. As the problem is convex with respect to each block of variables $U_k,V_k$ and $X_{[k]}$ while the other two are fixed, we have the following updates:
\bea
U^{l+1}_{k} &=& X^{l}_{[k]}(V^{l}_{k})^{T}(V^{l}_{k}(V^{l}_{k})^{T})^{\dagger},\label{U1}\\
V^{l+1}_{k} &=&( (U^{l+1}_{k})^{T}U^{l+1}_{k})^{\dagger}(U^{l+1}_{k})^{T})X^{l}_{[k]}\\
X^{l+1}_{[k]}&=& U^{l+1}_{k}V^{l+1}_{k},
\label{X1}
\eea
where \textquotedblleft$^\dagger$\textquotedblright denotes the Moore-Penrose pseudoinverse. It was shown in \cite{Xu} that, we can replace (\ref{U1}) by the following one:
\bea
U^{l+1}_{k} &=& X^{l}_{[k]}(V^{l}_{k})^{T},\label{U2}
\eea
to avoid computing the Moore-Penrose pseudoinverse $(V^{l}_{k}(V^{l}_{k})^{T})^{\dagger}$. The rationale behind this is that we only need the product $U^{l+1}_{k}V^{l+1}_{k}$ to compute $X^{l+1}_{[k]}$ as in (\ref{X1}) that is the same when either (\ref{U1}) or (\ref{U2}) is used. After updating $U^{l+1}_k,V^{l+1}_k$ and $X^{l+1}_{[k]}$ for all $k=1,2,\ldots,N-1$, we compute elements of the tensor $\mc{X}^{l+1}$ as follows:
\bea
x^{l+1}_{i_1\cdots i_N}= \left\{\begin{array}{ll}\Big(\sum_{k=1}^{N-1} \al_k \text{fold}(X^{l+1}_{[k]})\Big)_{i_1\cdots i_N}&({i_1\cdots i_N})\notin\Omega\\
	t_{i_1\cdots i_N}&({i_1\cdots i_N})\in\Omega
\end{array}\right.
\label{UpdateX2}
\eea
Let us name the algorithm as TMac-TT which stands for "tensor completion by parallel matrix factorization in the concept of tensor train" and its pseudo-code is summarized in Algorithm \ref{Algorithm 2}. Again, the Algorithm \ref{Algorithm 2} can be applied for the square model \cite{Mu_2014} by choosing the weights such that $\al_{k}=1$ if $k=\text{round}(N/2)$ otherwise $\al_{k}=0$, and we call it as TMac-Square.
\begin{table}[!thb]
	\centering
	\caption{TMac-TT}
	\label{Algorithm 2}	
	\begin{tabular}{*2l} % centered columns (4 columns)
		\hline
		~&~\\
		\tb{Input:} The observed data $\mc{T}\in\mathbb{R}^{I_1\times I_2\cdots\times I_{N}}$, index set $\Omega$.&\\
		\tb{Parameters:} $\al_i,r_i, i=1,\ldots,N-1$.&\\							
		~&~\\	
		\hline		
		~&~\\	
		1:~~\tb{Initialization:} $U^{0}, V^{0}, \mc{X}^{0}$, with $\mc{X}^{0}_{\Omega} = \mc{T}_{\Omega}$, $l=0$.& \\
		\tb{While not converged do:}&\\
		2:~~\tb{for} $k = 1$ \tb{to} $N-1$ \tb{do}&\\
		3:~~~~~~~Unfold the tensor $\mc{X}^{l}$ to get $X^{l}_{[k]}$&\\
		4:~~~~~~~$U^{l+1}_{i} = X^{l}_{[k]}(V^{l}_{k})^{T}$&\\
		5:~~~~~~~$V^{l+1}_{k} = ((U^{l+1}_{k})^{T}U^{l+1}_{k})^{\dagger}(U^{l+1}_{k})^{T}X^{l}_{[k]}$&\\	
		6:~~~~~~~$X^{l+1}_{[k]} = U^{l+1}_{k}V^{l+1}_{k}$&\\		
		7:~~\tb{end}&\\	
		8:~~Update the tensor $\mc{X}^{l+1}$ using (\ref{UpdateX2})\\
		\tb{End while}&\\
		\hline
		~&~\\
		\tb{Output:} The recovered tensor $\mc{X}$ as an approximation of $\mc{T}$&
	\end{tabular}
\end{table}
\subsection{Computational complexity of algorithms}
We analyze the computational complexity of algorithms applied to complete a tensor $\mc{X}\in\mathbb{R}^{I_1\times I_2\times\cdots\times I_N}$ in the Table \ref{Cost} where we assume that $I_1=I_2=\cdots= I_{N}=I$, the Tucker rank and TT rank are the same $r_1=r_2=\cdots= r_{N}=r$. 
\captionsetup[table]{name={\bf Table}}
\setcounter{table}{0}
\begin{table}[!h]
	\caption{Computational complexity of algorithms for one iteration.}
	\label{Cost}
	\centering % used for centering table
	\begin{tabular}{l l}
		Algorithm & Computational complexity\\
		\hline
		~\\
		SiLRTC  & $O(NI^{N+1})$\\
		SiLRTC-TT & $O(I^{3N/2}+I^{3N/2-1})$\\
		TMac & $O(3NI^Nr)$\\
		TMac-TT & $O(3(N-1)I^Nr)$\\
		\hline
	\end{tabular}
\end{table}
\section{Tensor augmentation \label{sec3}}
In this section, we introduce the \emph{ket augmentation} (KA) to represent a lower-order tensor by a higher-order one, i.e.
 to  cast an $N$th-order tensor
$\mc{T}\in\mathbb{R}^{I_1\times I_2\times\cdots\times I_N}$ into a $K$th-order tensor
$\tilde{\mc{T}}\in\mathbb{R}^{J_1\times J_2\times\cdots\times J_K}$, where $K\geq N$ and $\prod_{l=1}^{N}I_l=\prod_{l=1}^{K}J_l$. Higher-order representation of the tensor offers some important advantages. For instance, TT decomposition is more efficient for
the augmented tensor  because  the local structure of the data can be exploited effectively in terms of computational resources.
Actually, if the tensor is slightly correlated, its augmented tensor can be represented by a low-rank TT
\cite{Oseledets_2011,Lattore}.

The KA was originally introduced in \cite{Lattore} for casting a grayscale image into \emph{real ket state} of a Hilbert space, which is simply a higher-order tensor, using an appropriate block structured addressing. Here we generalize the KA scheme for third-order tensors $\mc{T}\in\mathbb{R}^{I_1\times I_2\times I_{3}}$ that represent color images, where $I_1\times I_2=2^n\times 2^n$ ($n\geq 1\in\mathbb{Z}$) is the number of pixels in the image and $I_3 = 3$ is the number of colors (red, green and blue).
Let us start with an initial block, labeled as $i_1$, of $2\times 2$ pixels corresponding to a single color $j$ (assume that the color is indexed by $j$ where $j=1,2,3$ corresponding to red, green and blue colors, respectively). This block can be represented as
\bea
\mc{T}_{[2^1\times 2^1\times j]} = \sum_{i_1=1}^{4}c_{i_1j}\tb{e}_{i_1},
\eea
where $c_{i_1j}$ is the pixel value corresponding to color $j$ and $\tb{e}_{i_1}$ is the orthonormal base which is defined as $\tb{e}_{1} = (1,0,0,0)$, $\tb{e}_{2} = (0,1,0,0)$, $\tb{e}_{3} = (0,0,1,0)$ and $\tb{e}_{4} = (0,0,0,1)$. The value $i_1=1,2,3$ and $4$ can be considered as labeling the up-left, up-right, down-left and down-right pixels, respectively. For all three colors, we have three blocks which are presented by
\bea
\mc{T}_{[2^1\times 2^1\times 3]} = \sum_{i_1=1}^{4}\sum_{j=1}^{3}c_{i_1j}\tb{e}_{i_1}\otimes \tb{u}_{j},
\label{ket1}
\eea
where $\tb{u}_{j}$ is also an orthonormal base which is defined as $\tb{u}_{1} = (1,0,0)$, $\tb{u}_{2} = (0,1,0)$, $\tb{u}_{3} = (0,0,1)$.
\begin{figure}[htpb]
	\centering
	\includegraphics[width=\columnwidth]{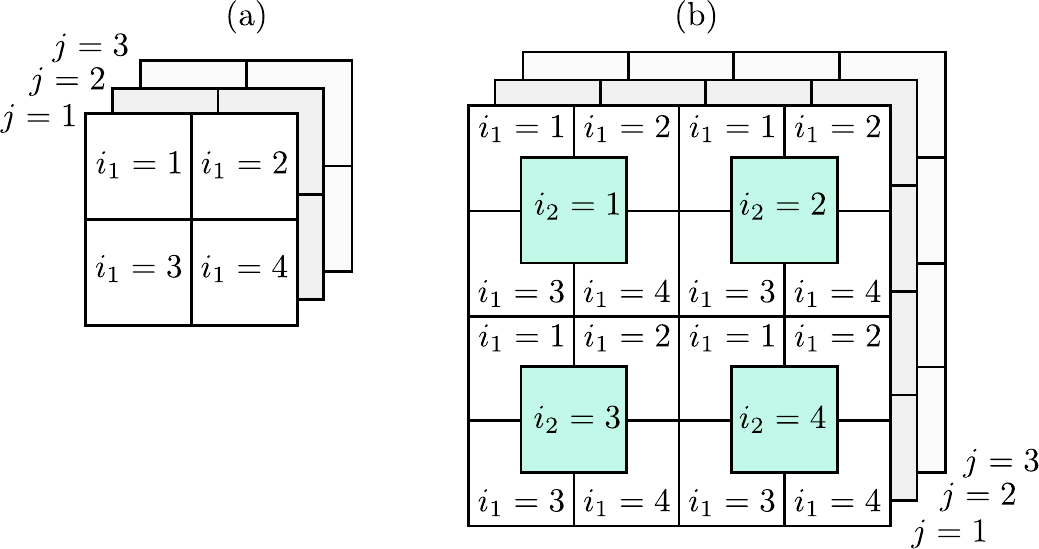}\\
	\caption{A structured block addressing procedure to cast an image into a higher-order tensor. (a) Example for an image of size $2\times 2\times 3$ represented by (\ref{ket1}). (b) Illustration for an image of size $2^2\times 2^2\times 3$ represented by (\ref{ket2}).}
	\label{fig10}
\end{figure}
We now consider a larger block labeled as $i_2$ make up of four inner sub-blocks for each color $j$ as shown in Fig.~\ref{fig10}. In total, the new block is represented by
\bea
\mc{T}_{[2^2\times 2^2\times 3]} = \sum_{i_2=1}^{4}\sum_{i_1=1}^{4}\sum_{j=1}^{3}c_{i_2i_1j}\tb{e}_{i_2}\otimes\tb{e}_{i_1}\otimes \tb{u}_{j}.
\label{ket2}
\eea
Generally, this block structure can be extended to a size of $2^n\times 2^n\times 3$ after several steps until it can present all the values of pixels in the image. Finally, the image can be cast into an $(n+1)$th-order tensor $\mc{C}\in\mb{R}^{4\times 4\times\cdots\times 4\times 3}$ containing all the pixel values as follows,
\bea
\mc{T}_{[2^n\times 2^n\times 3]} = \sum_{i_n,\ldots,i_1=1}^{4}\sum_{j=1}^{3}c_{i_n\cdots i_1j}\tb{e}_{i_n}\otimes\cdots\otimes\tb{e}_{i_1}\otimes \tb{u}_{j}.
\eea

When the image is represented by a real ket state, its entanglement entropy can reflect the correlation between individual pixels as due to their relative positions in the image. Besides, this presentation is suitable for the image processing as it not only preserves the pixels values of the image but also rearrange them in a higher-order tensor such that the richness of textures in the image can be studied via the correlation between modes of the tensor \cite{Lattore}.
\section{Simulations \label{sec4}}
We apply the proposed algorithms for completing both synthetic data and color images. Simulation results are compared with commonly used algorithms, i.e. SiLRTC \cite{Ji2013}, SiLRTC-Square \cite{Mu_2014}, TMac \cite{Xu} and TMac-Square. To measure performance of a LRTC algorithm we compute the relative square error (RSE) between the approximately recovered tensor $\mc{X}$ and the original one $\mc{T}$, which is defined as,
\bea
RSE = ||\mc{X}-\mc{T}||_F/||\mc{T}||_F.
\label{RSE}
\eea
It is hard to choose the parameters in the models so that optimal solutions can be achieved. In experiments, we simply choose the weights $\al_k$ for our proposed algorithms as follows:
\bea
\al_k &=& \frac{\delta_{k}}{\sum_{k=1}^{N-1}\delta_{k}}~~\text{with}~~\delta_{k} = \min(\prod_{l=1}^{k}I_l,\prod_{l=k+1}^{N}I_l),
\eea
where $k=1,\ldots,N-1$. The positive parameters are chosen by $\bet_k = f\al_k$, where $f$ is empirically chosen from one of the following values in $[0.01, 0.05, 0.1, 0.5, 1]$ in such a way that the algorithm performs the best. For the algorithms used to compare with ours, i.e. SiLRTC and TMac the weights are chosen as follows:
\bea
\al_k &=& \frac{I_{k}}{\sum_{k=1}^{N}I_{k}},
\eea
where  $k=1,\ldots,N$. The positive parameters are chosen such that $\bet_k = f\al_k$, where $f$ is empirically chosen from one of the following values in $[0.01, 0.05, 0.1, 0.5, 1]$ which gives the best performance. The convergence criterion of our proposed algorithms is defined by computing the relative error of the tensor $\mc{X}$ between two successive iterations as follows:
\bea
\epsilon = \frac{||\mc{X}^{l+1} - \mc{X}^{l}||_{F}}{||\mc{T}||_{F}}\leq tol,
\eea
where we set $tol=10^{-4}$ and the maximum number of iterations is $maxiter = 1000$. 

In what follows, we perform simulations for algorithms with respect to different missing ratios ($mr$) defined as,
\bea
mr = \frac{p}{\prod_{k=1}^{N}I_k},
\eea
where $p$ is the number of missing entries which are often chosen randomly from the tensor $\mc{T}$ based on a uniform distribution. These simulations are implemented under a Matlab environment using the FEIT cluster from the University of Technology Sydney.
\subsection{Synthetic data completion}
We firstly perform the simulation on two different types of low-rank tensors which are generated synthetically in such a way that the Tucker and TT rank are known in advance.
\subsubsection{Completion of low TT rank tensor}
The $N$th-order tensors $\mc{T}\in\mathbb{R}^{I_1\times I_2\cdots\times I_{N}}$ of TT rank $(r_1,r_2,\ldots, r_{N-1})$ are generated such that its elements is represented by a TT format \cite{Oseledets_2011}. Specifically, its elements is $t_{i_{1}i_{2}\ldots i_{N}} = A^{[1]}_{i_1}A^{[2]}_{i_2}\cdots A^{[N]}_{i_N}$, where $A^{[1]}\in\mathbb{R}^{I_1\times r_{1}}$, $A^{[N]}\in\mathbb{R}^{r_{N}\times I_N}$ and $\mc{A}^{[k]}\in\mathbb{R}^{r_{k-1}\times I_k\times r_{k}}$ with $k = 2,\ldots,N-1$ are generated randomly with respect to the standard Gaussian distribution $\mc{N}(0,1)$. For simplicity, in this paper we set all components of the TT rank the same and so does the dimension of each mode, i.e. $r_1=r_2=\cdots= r_{N-1}=r$ and $I_1=I_2=\cdots= I_{N}=I$.

\begin{figure}[htpb]
	\centering
	\includegraphics[width=\columnwidth]{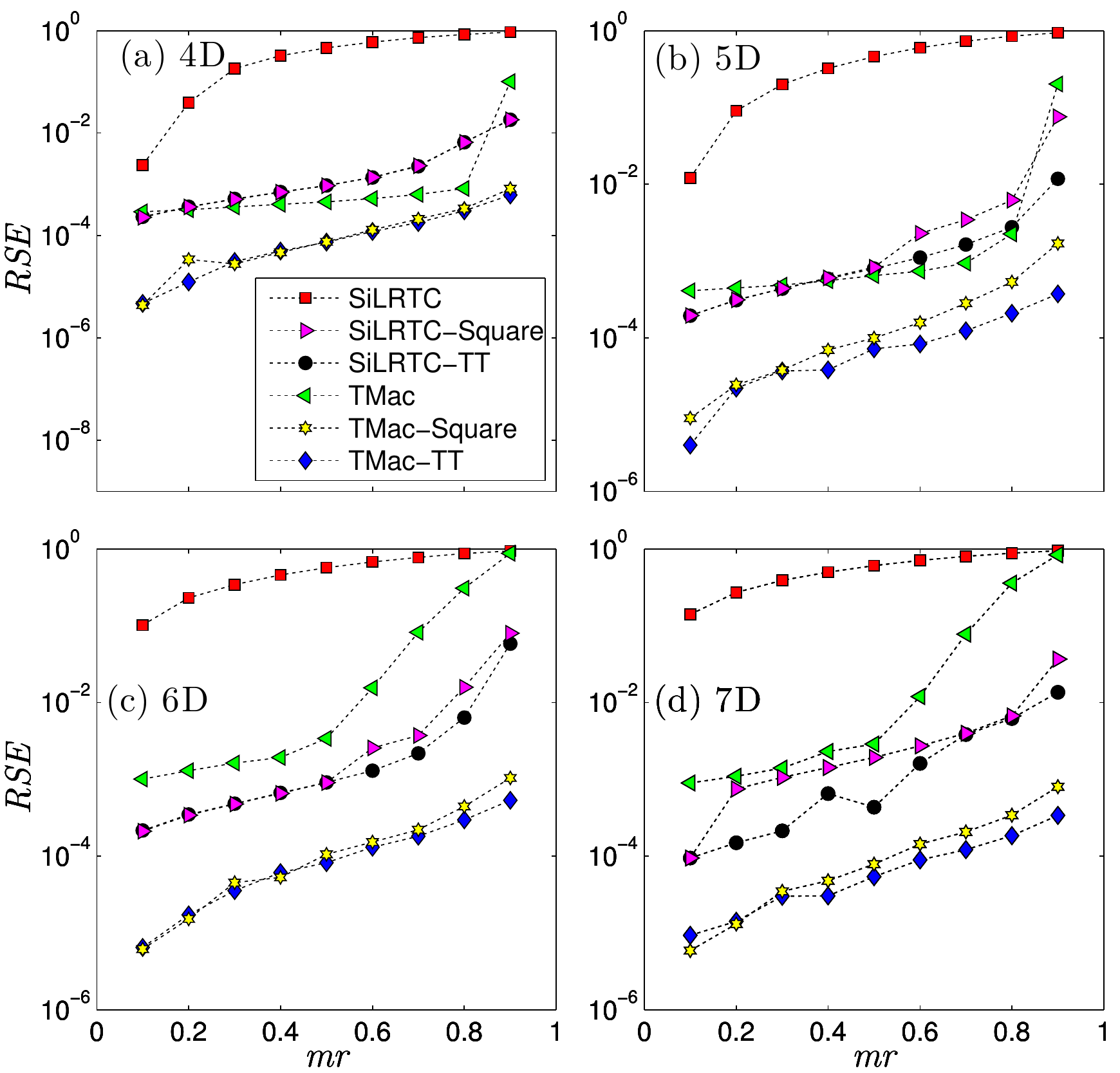}\\
	\caption{The RSE comparison when applying different LRTC algorithms to synthetic random tensors of low TT rank. Simulation results are shown for different tensor dimensions, 4D, 5D, 6D and 7D.}
	\label{fig1}
\end{figure}

The plots of RSE with respect to $mr$ are shown in the Figure.~\ref{fig1} for tensors of different sizes, $40\times 40\times 40\times 40$ (4D), $20\times 20\times 20\times 20\times 20$ (5D), $10\times 10\times 10\times 10\times 10\times 10$ (6D) and $10\times 10\times 10\times 10\times 10\times 10\times 10$ (7D) and the corresponding TT rank tuples are $(10, 10, 10)$ (4D), $(5, 5, 5,5)$ (5D), $(4,4,4,4,4)$ (6D) and $(4,4,4,4,4,4)$ (7D). From the plots we can see that TMac-TT shows best performance in most cases. Especially, TMac-TT can recover the tensor successfully despite the high missing ratios. Particularly, in most cases with high missing ratios, e.g. $mr = 0.9$, it can recover the tensor with $RSE\approx 10^{-4}$. More importantly, the proposed algorithms SiLRTC-TT and TMac-TT often performs better than their corresponding counterparts, i.e. SiLRTC and TMac in most cases.

\begin{figure}[htpb]
	\centering
	\includegraphics[width=\columnwidth]{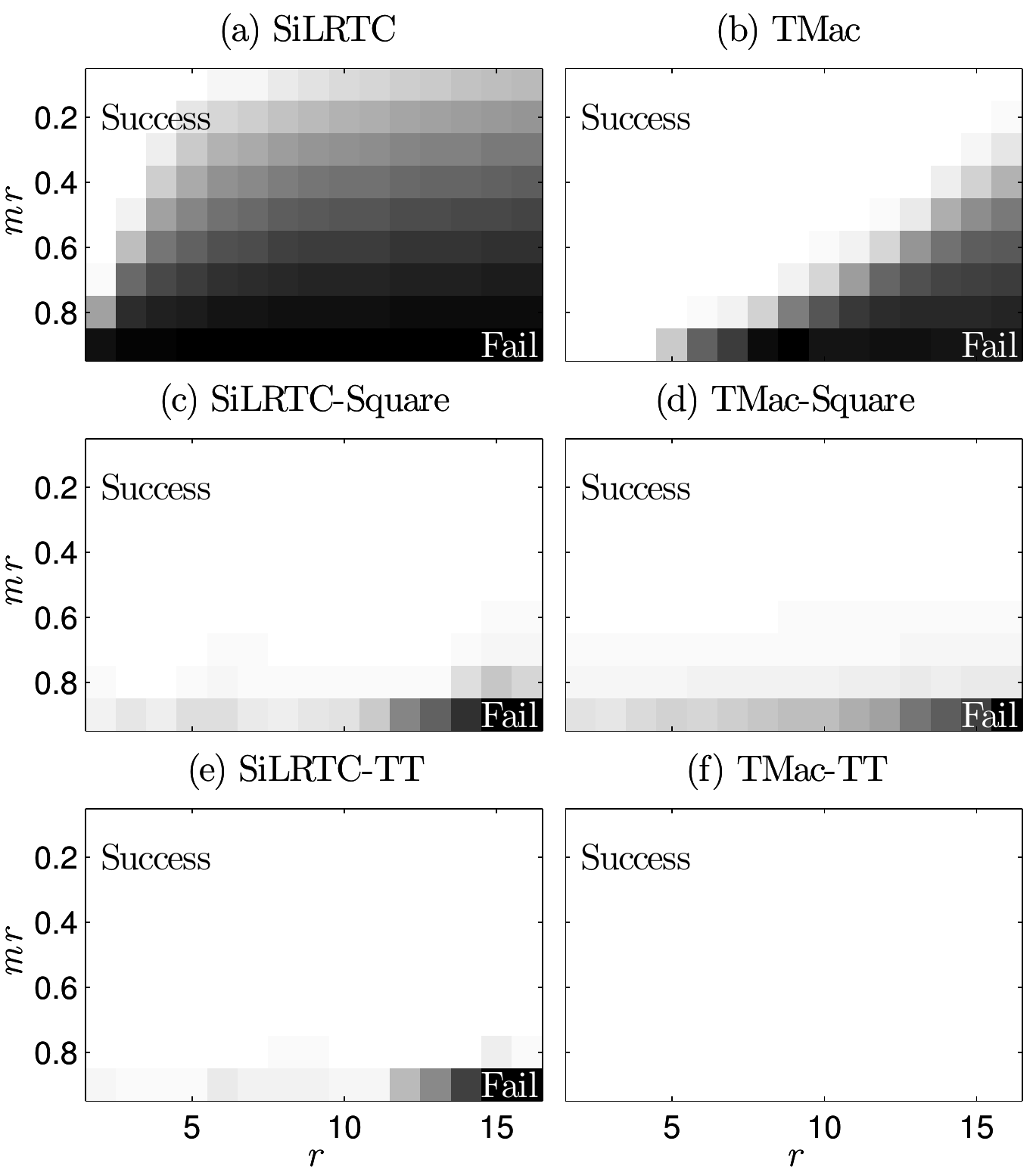}\\
	\caption{Phase diagrams for low TT rank tensor completion when applying different algorithms to a 5D tensor.}
	\label{fig2}
\end{figure}

For a better comparison on the performance of different LRTC algorithms, we present the phase diagrams using the grayscale color to estimate how successfully a tensor can be recovered for a range of different TT rank and missing ratios. If $RSE\leq\epsilon$ where $\epsilon$ is a small threshold, we say that the tensor is recovered successfully and is represented by a white block in the phase diagram. Otherwise, if $RSE>\epsilon$, the tensor is recovered partially with a relative error and the block color is gray. Especially the recovery is completely failed if $RSE=1$. Concretely, we show in Fig.~\ref{fig2} the phase diagrams for different algorithms applied to complete a 5D tensor of size $20\times 20\times 20\times 20\times 20$ where the TT rank $r$ varies from 2 to 16 and $\epsilon=10^{-2}$. We can see that our LRTC algorithms outperform the others. Especially, TMac-TT always recovers successfully the tensor with any TT rank and missing ratio.
\subsubsection{Completion of low Tucker rank tensor}
\begin{figure}[htpb]
	\centering
	\includegraphics[width=\columnwidth]{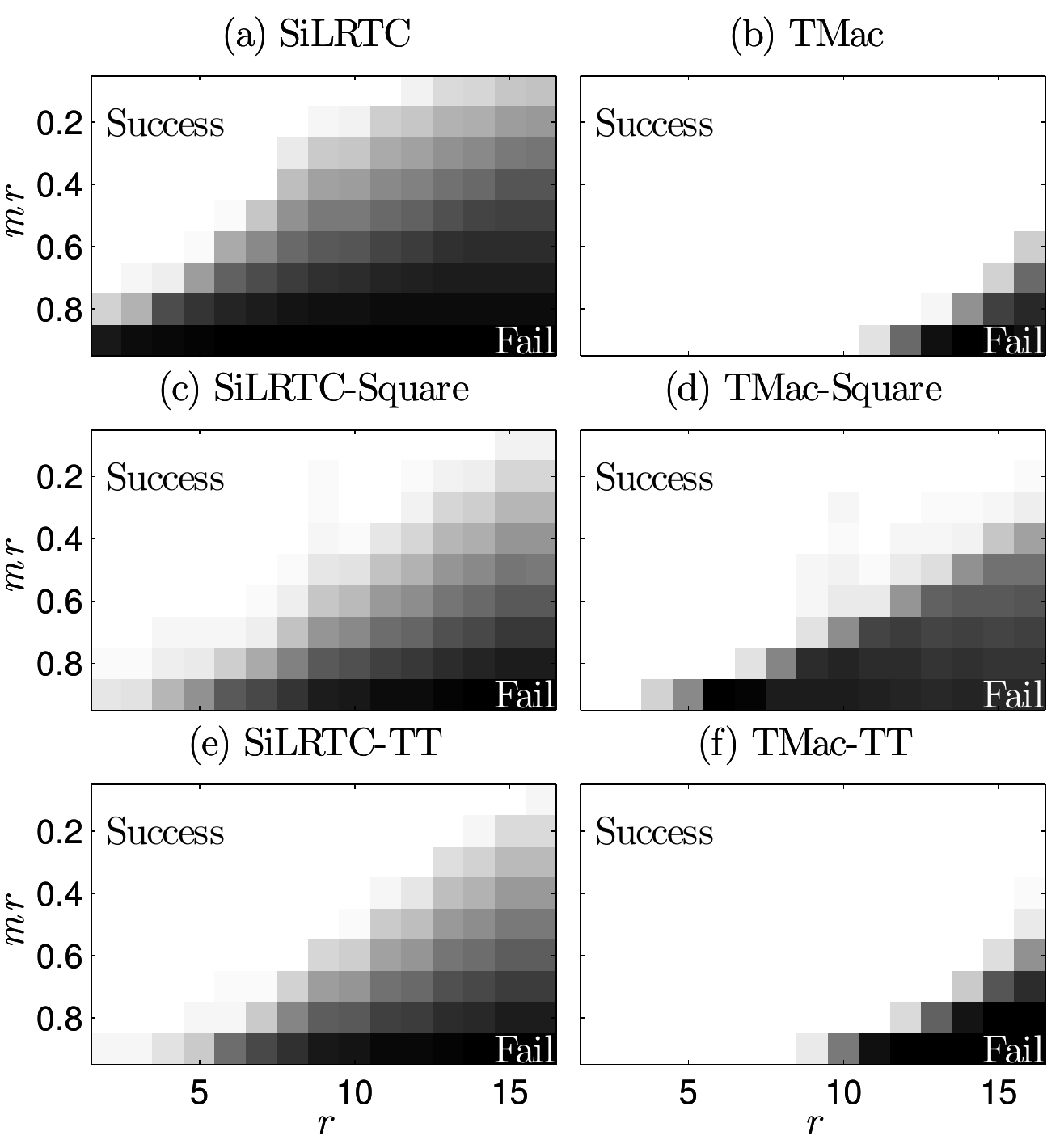}\\
	\caption{Phase diagrams for low Tucker rank tensor completion when applying different algorithms to a 5D tensor.}
	\label{fig4}
\end{figure}

Let us now apply our proposed algorithms to synthetic random tensors of low Tucker rank. The $N$th-order tensor $\mc{T}\in\mathbb{R}^{I_1\times I_2\cdots\times I_{N}}$ of Tucker rank $(r_1,r_2,\ldots, r_N)$ is constructed by $\mc{T} = \mc{G}\times_1A^{(1)}\times_2A^{(2)}\cdots\times_NA^{(N)}$, where the core tensor $\mc{G}\in\mathbb{R}^{r_1\times r_2\cdots\times r_{N}}$ and the factor matrices $A^{(k)}\in\mathbb{R}^{r_k\times I_{k}},k=1,\ldots,N$ are generated randomly by using the standard Gaussian distribution $\mc{N}(0,1)$. Here, we choose $r_1=r_2=\cdots= r_{N}=r$ and $I_1=I_2=\cdots= I_{N}=I$ for simplicity. To compare the performance between the algorithms, we show in the Fig.~\ref{fig4} the phase diagrams for different algorithms applied to complete a 5D tensor of size $20\times 20\times 20\times 20\times 20$ where the Tucker rank $r$ varies from 2 to 16 and $\epsilon=10^{-2}$. We can see that both TMac and TMac-TT perform much better than the others and. Besides, SiLRTC-TT shows better performance when compared to SiLRTC and SiLRTC-Square. Similarly, TMac-TT is better than its particular case TMac-Square.

In summary, we can see that although the tensors are generated synthetically to have low Tucker ranks, the proposed algorithms are still capable of producing results which are as good as the ones obtained by the Tucker-based algorithms. In order to have a better comparison between algorithms, we show results of applying them to the real world data such as color images where the ranks of the tensors are not known in advance in the next subsection.
\subsection{Image completion}
A set of color images, namely "Peppers", "Lenna" and "House" are employed to test the algorithms with different missing ratios. All the images are initially represented by third-order tensors which have same sizes of $256\times 256\times 3$.
\begin{figure}[htpb]
	\centering
	\includegraphics[width=\columnwidth]{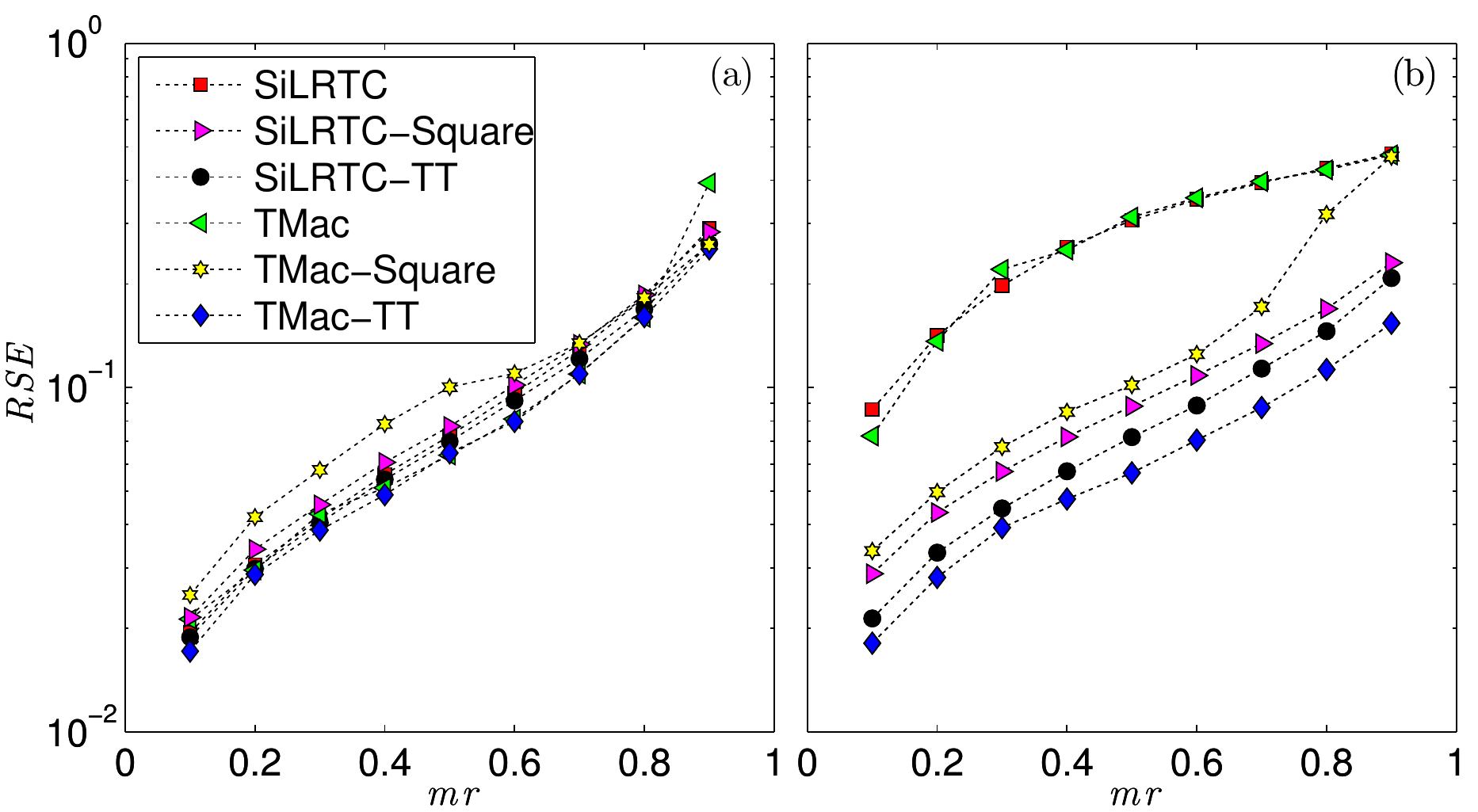}\\
	\caption{Performance comparison between different tensor completion algorithms based on the RSE vs the missing rate when applied to the Peppers image. (a) Original tensor (no order augmentation). (b) Augmented tensor using KA scheme.}
	\label{fig5}
\end{figure}
Note that when completing the third-order tensors, we do not expect that our proposed methods prevail against the conventional ones due to the fact that the TT rank of the tensor is a special case of the Tucker rank. Thus, performance of the algorithms should be mutually comparable. However, for the purpose of comparing the performance between different algorithms for real data (images) represented in terms of higher-order tensors, we apply tensor augmentation scheme KA mentioned above to reshape third-order tensors to  higher-order ones without changing the number of entries in the tensor. Specifically, we start our simulation by casting a third-order tensor $\mc{T}\in\mathbb{R}^{256\times 256\times 3}$ into a ninth-order $\tilde{\mc{T}}\in\mathbb{R}^{4\times  4\times 4\times 4\times 4\times  4\times 4\times 4\times 3}$ and then applying the tensor completion algorithms to impute its missing entries. We perform the simulation for the Peppers and Lenna images where missing entries of each image are chosen randomly according to a uniform distribution, the missing ratio $mr$ varies from 0.1 to 0.9. In Fig.~\ref{fig5}, we compare performance of algorithms on completing the Peppers image. We can see that, when the image is represented by a third-order tensor, the performance of the algorithms are comparable (The TMac-TT is actually slightly better than the others in most of the missing ratios). However, for the case of the ninth-order tensors, the performance of the algorithms are rigorously distinguished. Specifically, our proposed algorithms (especially TMac-TT) prevail against the others. We also illustrate the recovered images for $mr=0.7$ in Fig.~\ref{fig6}. This shows that our proposed algorithms give really good results in the case of augmented tensors, meanwhile the compared algorithms seem to be inefficient. Furthermore, using the KA scheme to increase the tensor order, SiLRTC-TT and TMac-TT significantly improve the accuracy when compared to the cases without augmentation. More precisely, TMac-TT gives the best results $RSE\approx0.088$ with respect to the case of using KA scheme. Same experiment is performed on the Lenna image and recovery results are shown in Fig.~\ref{fig7} and Fig.~\ref{fig8}. The results also show that TMac-TT gives the best results for the augmented tensor using the KA scheme. 
\begin{figure}[htpb]
	\centering
	\includegraphics[width = \columnwidth]{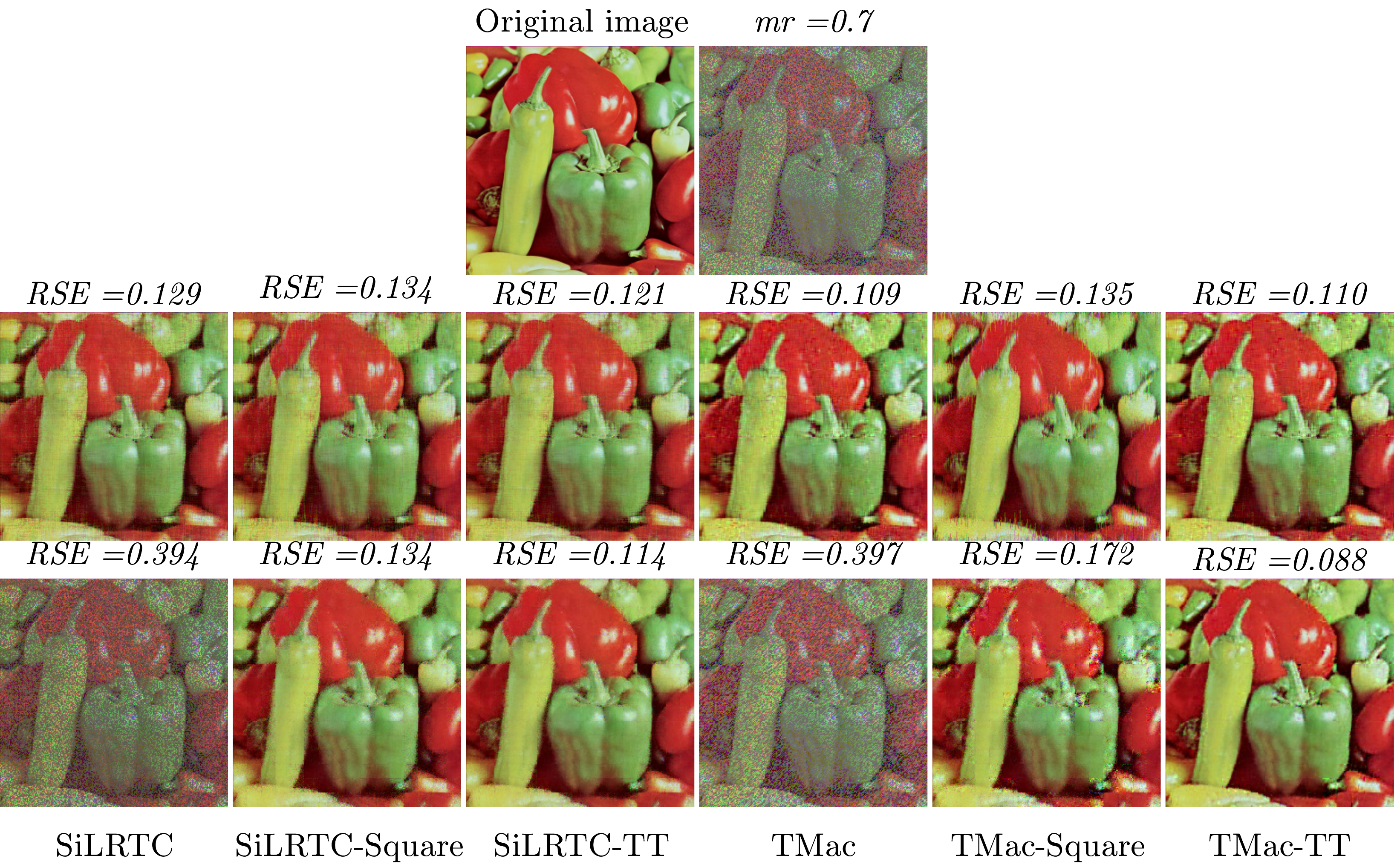}\\
	\caption{Recover the Peppers image with $70\%$ of missing entries using different algorithms. Top row from left to right: the original image and its copy with $70\%$ of missing entries. Second and third rows represent the recovery results of third-order (no order augmentation) and ninth-order tensors (KA augmentation),  using different algorithms: SiLRTC, SiLRTC-Square, SiLRTC-TT, TMac, TMac-Square and TMac-TT from the left to the right, respectively.}
	\label{fig6}
\end{figure}
\begin{figure}[htpb]
	\centering
	\includegraphics[width=\columnwidth]{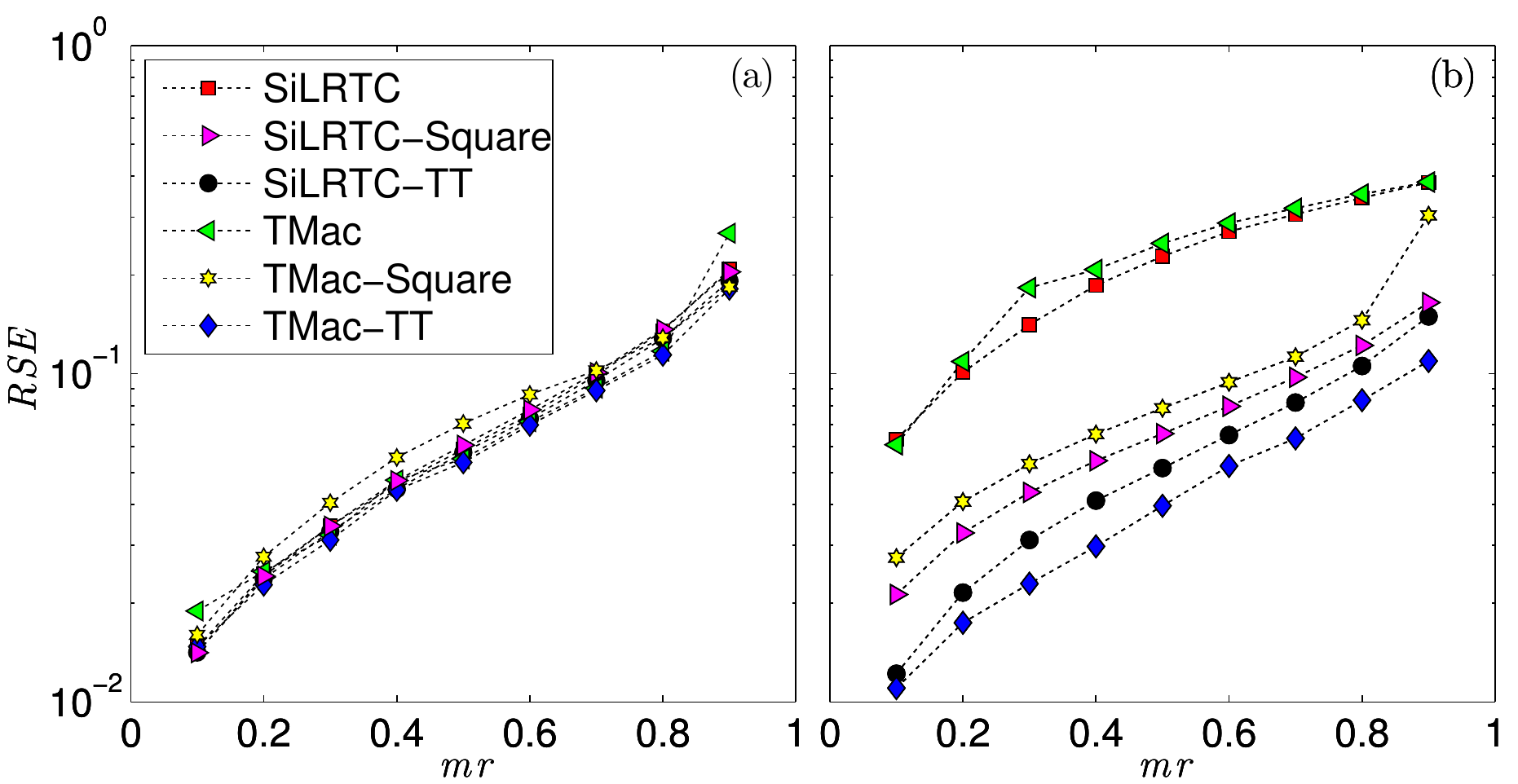}\\
	\caption{Performance comparison between different tensor completion algorithms based on the RSE vs the missing rate when applied to the Lenna image. (a) Original tensor (no order augmentation). (b) Augmented tensor using KA scheme.}
	\label{fig7}
\end{figure}
\begin{figure}[htpb]
	\centering
	\includegraphics[width = \columnwidth]{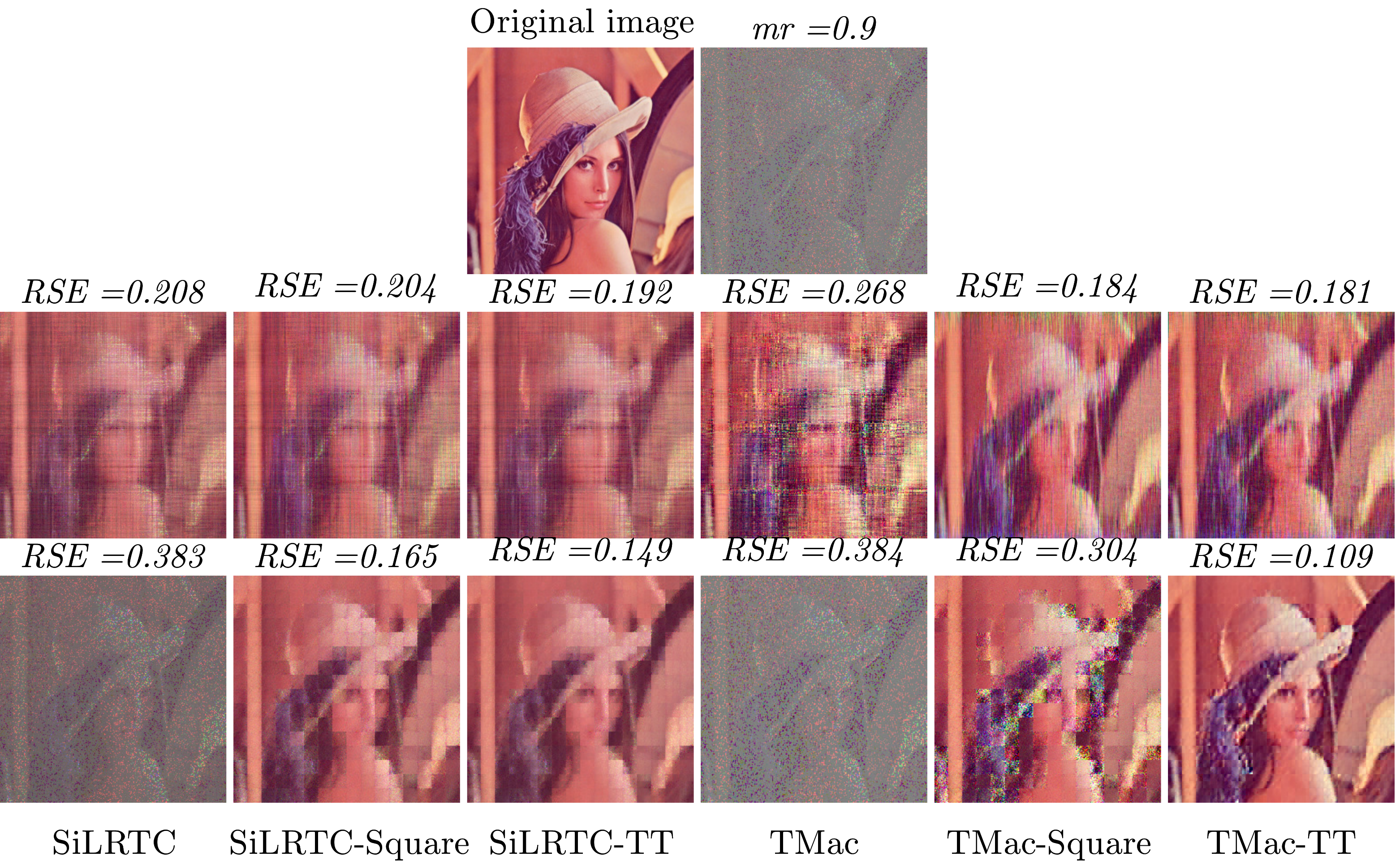}\\
	\caption{Recover the Lenna image with $90\%$ of missing entries using different algorithms. Top row from left to right: the original image and its copy with $90\%$ of missing entries. Second and third rows represent the recovery results of third-order (no order augmentation) and ninth-order tensors (KA augmentation),  using different algorithms: SiLRTC, SiLRTC-Square, SiLRTC-TT, TMac, TMac-Square and TMac-TT from the left to the right, respectively.}
	\label{fig8}
\end{figure}

\begin{figure}[htpb]
	\centering
	\includegraphics[width = \columnwidth]{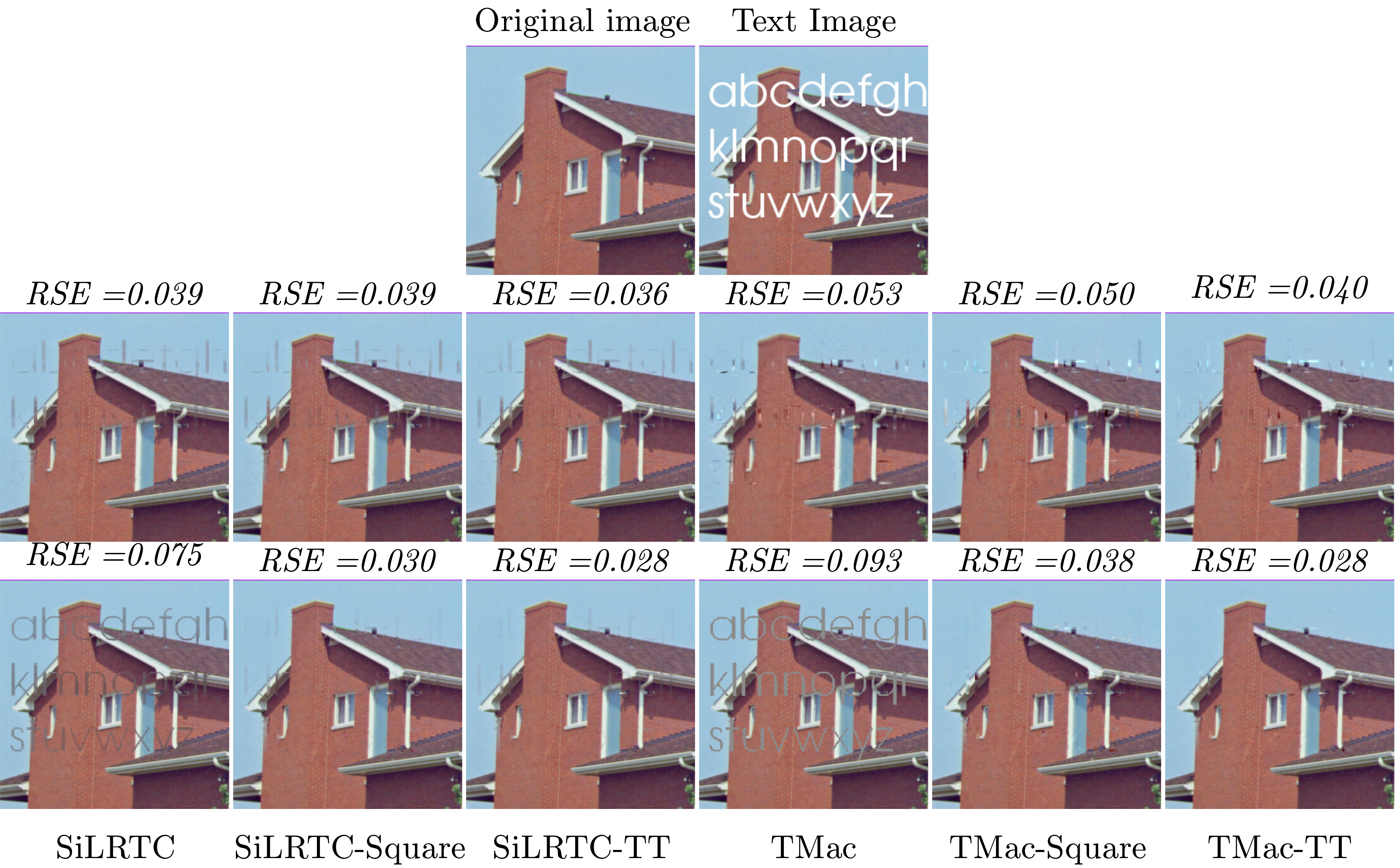}\\
	\caption{Recover the House image with missing entries described by the white letters using different algorithms. Top row from left to right: the original image and its copy with white letters. Second and third rows represent the recovery results of third-order (no order augmentation) and ninth-order tensors (KA augmentation), using different algorithms: SiLRTC, SiLRTC-Square, SiLRTC-TT, TMac, TMac-Square and TMac-TT from the left to the right, respectively.}
	\label{fig9}
\end{figure}

We perform the same above experiment on the House image, however, the missing entries are now chosen as the white text, and hence the missing rate is fixed. The result is shown in Fig.~\ref{fig9}. In the cases of tensor augmentation, the conventional algorithms SiLRTC and TMac do not perform well meanwhile our proposed algorithms do. Using the KA scheme, better results can be achieved by employing our algorithms when compared to the case without using the augmentation schemes.

To sum up, we see that the TT-based algorithms outperforms the Tucker-based ones when applying to the images represented by ninth-order tensors $\tilde{\mc{T}}\in\mathbb{R}^{4\times  4\times 4\times 4\times 4\times  4\times 4\times 4\times 3}$. This is because the components of TT rank can approximately vary in a broad range of values (the maximum value it can reach is 256) to capture the global information of the images. On the contrary, the components of Tucker rank can have value up to $4$ due to the mode-$k$ matricization $X_{{k}}\in\mb{R}^{4\times 49152}$. Consequently, the Tucker-based algorithms are not reliable due to the naturally small Tucker rank.

\section{Conclusion \label{sec5}}
We have proposed efficient LRTC algorithms based on the concept of the TT decomposition. The SiLRTC-TT algorithm is applied to minimize the TT rank of the tensor by solving the TT nuclear norm optimization. Meanwhile, TMac-TT is based on the multilinear matrix factorization model to minimize the TT-rank. The latter is more computationally efficient due to the fact that it does not need the SVD which is different from the former. The proposed algorithms are employed to simulate with both synthetic and real world data represented by higher-order tensors and their performance are compared with their replicates, which are formulated in terms of Tucker rank. For synthetic data, on the one hand our algorithms prevail the other when the tensors have low TT rank. On the other hand, their performance are comparable in case of low Tucker rank tensors. Therefore, the TT-based algorithms are quite promising and reliable when applying to real world data. To validate this, we apply the algorithms to study the image completion problem. Benchmark results show that when applied to original tensors without order augmentation, all algorithms are comparable to each other. However, in the case of augmented tensors, our proposed algorithms not only outperform the others but also provide better recovery results when compared to the case without tensor order augmentation.

Although the proposed algorithms can potentially be applied to complete tensors with a wide range of low tensor ranks, i.e. Tucker rank or TT rank, the optimal parameters such as weights and TT rank cannot be chosen automatically rather than empirically. We plan to further improve the algorithms by developing a scheme to adaptively choose these parameters. Besides, their applications
to data compression, text mining, image classification and video indexing are under our interest.

\bibliographystyle{IEEEtran}
\bibliography{TC}
\end{document}